\documentclass[a4paper,11pt]{amsart} 
\usepackage{amsmath, amsthm, amssymb}

\usepackage[inner=3cm,outer=3cm,top=2.8cm,
bottom=2.8cm,marginparwidth=90pt]{geometry} 
\usepackage{times}

\newtheorem{theorem}{Theorem}%[chapter]

\newtheorem{lemma}[theorem]{Lemma}

\theoremstyle{definition}
\newtheorem*{remark}{{Remark}}
\newtheorem*{example}{Example}

\newcommand{\bigR}{{\mathbb R}}

\DeclareMathOperator{\trace}{trace}

\newcommand{\pd}[2]{\dfrac{\partial#1}{\partial#2}}

\newcommand{\ktwo}{k}
\newcommand{\kone}{k'}

\newcommand{\eff}{F} % the area integrand- either write as \bar F, or just F 
\newcommand{\normal}{\nu} %notation for the normal--- either \omega or \nu...

\title{Interior gradient estimates for anisotropic mean curvature flow}
\author{Julie Clutterbuck}
\def\fuhome{@math.fu-berlin.de}
\address{Freie Universit\"at Berlin, Arnimallee 2-6,
  14195 Berlin, Germany}
\email{Julie.Clutterbuck\fuhome}


\begin{document}

%%%% ABSTRACT  %%%%%%%%%%%%
\begin{abstract}
In this paper we consider the evolution of a graph-like hypersurface by anisotropic mean curvature flow, under some restrictions on the anisotropic area integrand.  We find interior estimates (in both time and space) on the gradient of such hypersurfaces, depending only on the height of the graph and the anisotropic area integrand.
\end{abstract}

\maketitle

%%%%% INTRODUCTION  %%%%%
\section{Introduction}

Consider the evolution of a hypersurface by its mean curvature:
 \begin{equation} \frac{d }{dt}
\mathbf{x}(y,t)=\mathbf{H}(y,t), \qquad y\in \mathbf{M}, \label{mean curvature flow} \end{equation}
where $\mathbf{x}:\mathbf{M}^n\times[0,T]\rightarrow \bigR^{n+1}$ is the immersion of a manifold $\mathbf{M}$ at each time $t$ and
$\mathbf{H}$ is the mean curvature vector.

$\mathbf{M}$ can be written as a graph when a fixed vector $\omega\in\bigR^{n+1}$ can be found so that for a choice of
unit normal $\normal$, $ %\begin{equation*} 
\langle \normal,\omega \rangle >0  %\end{equation*} 
$ everywhere.  %Equivalently, $\langle \nu,\omega \rangle ^{-1}$ is bounded above.   
Given the image $\mathbf{x}(y,t)$ of a point $y\in \mathbf{M}$,  the
{\emph{height}} of $\mathbf{M}$ above the hyperplane defined by $\omega$ is denoted by $
u=\langle \mathbf{x}, \omega \rangle$,
and the \emph{gradient function} is given by 
$ v=\langle \nu, \omega \rangle^{-1}=\sqrt{1+|Du|^2}$.  

In \cite{eh:mean},  Ecker and Huisken established that when initial data is given by an entire Lipschitz graph with a linear growth bound, there is a smooth solution to \eqref{mean curvature flow} for all times.  An important step in this proof was showing that the solutions remain graphs:  this was done by showing that $v$ is bounded above, with a constant depending on the initial Lipschitz bound.
In \cite{eh:interior}, it was established that the Lipschitz bound need only be local.    The estimates in this paper are intended in the spirit of the local gradient estimates of Section 2 of the latter paper.

Such gradient estimates may be found even if  the initial data is not Lipschitz. 
In \cite{evans-spruck-III}, Evans and Spruck showed that under mean curvature flow, surfaces that may initially be written locally as a continuous graph, become smooth for $t>0$.  The level set method was also used by Barles, Biton and Ley \cite{barles:2002} to find similar gradient estimates for more a general class of equations.  More recently, in \cite{colding:2004}, Colding and Minicozzi found an explicit local estimate in the form  
\begin{equation*}
|Du(x,t)|\le \exp\left\lbrace{c\left(1+t^{-1/2}\|u\|_\infty\right)^2}\right\rbrace,
\end{equation*}
for solutions over a ball $B_R(x)$, $R\ge C\sqrt{t}$.  Here, the constant $c$  depends only the dimension.    This estimate does not depend on an initial gradient estimate.

In this paper we find an analogous result for anisotropic mean curvature flow, under some restrictions on the anisotropy.   

Such estimates are an important step in finding existence results for a variety of boundary value problems with non-smooth initial data, as in \cite{andrews}. 

In this paper, we follow the exposition of Andrews in \cite{andrews:anisotropic}, in particular Section 8, in which the evolution equation for graph-like surfaces is derived.

  We consider surfaces $\mathbf{M}$ with local embeddings $x=y^ie_i+u(y^1,\dots,y^n)e_0$, and normal $\normal=Du-\phi^0=\sum_{i=1}^n u_i\phi^i-\phi^0$, where $\lbrace \phi^0,\phi^1,\dots,\phi^n\rbrace$ and  $\lbrace e_0,e_1,\dots,e_n\rbrace$ are dual  bases for the cotangent space $V^*\cong\bigR^{n+1}$ and tangent space $V\cong\bigR^{n+1}$ respectively.   

The anisotropic mean curvature flow for such a surface is given by 
\begin{equation}
u_t=\left.{\eff}D^2{\eff}\right|_{Du-\phi^0} (\phi^i,\phi^j)D^2u(e_i,e_j)  \label{anisotropic curvature flows}
\end{equation}
(the homogeneous degree zero mobility function $m=m(\normal)$ of \cite{andrews:anisotropic} is here taken to be identically 1).    The anisotropic area integrand $\eff:V^*\rightarrow \bigR$ is a positive, convex  function that is $C^3$ on $V^*\backslash \lbrace 0\rbrace $, and homogeneous of degree one, so that $\eff (\lambda v)=\lambda\eff (v)$ for $v\in V^*$ and all scalars $\lambda>0$.  The level sets of $\eff $ are denoted by $\Sigma_\lambda:=\lbrace v \in V^*: \eff (v)=\lambda \rbrace$.  We impose the condition that the convex hull of each level set,  $\lbrace v  \in V^*: \eff (v )\le\lambda \rbrace$, must be uniformly convex.   

\begin{example} 
In the \emph{isotropic} case, $\eff(v)^2={\sum_{i=0}^n (v_i)^2}$ and the level sets of $\eff$ are spheres.  The coefficients are those of mean curvature: $\left.\eff D^2\eff\right|_{Du-\phi^0}(\phi^i,\phi^j)=\delta_{ij}-u_i u_j (1+|Du|^2)^{-1}$.  
\end{example}

We define two further conditions on $\eff$: 

% firstly, that third derivatives are small; secondly, that there is a symmetry in the distinguished direction $\phi^0$; and introduce a technical consequence of the symmetry condition.

\noindent{\textbf{Smallness of third derivatives condition: }}{\textit{ Define
\begin{equation*}
{Q_\normal}(p,q,r):= \eff^2(\normal) \left.D^3 \eff \right|_\normal(p,q,r). \label{defn of Q}
\end{equation*}
% , so that $\left.D\eff\right|_\normal(p)=0$ (and similarly for \dots)
This condition is satisfied when 
\begin{equation}
Q_\normal(p,q,r)\le C_1\left[\eff(\normal)^3\left.D^2\eff\right|_\normal(p,p)\left.D^2\eff\right|_\normal(q,q)\left.D^2\eff\right|_\normal(r,r)\right]^{1/2} \label{condition on D^3 F}
\end{equation}
for all $p,q,r$  tangent to the level set $\Sigma_{\eff (\normal)}$ at $\normal$. 
%tangent to the unit ball of $\eff$ at $\normal$, wh
Here $C_1$ is a positive constant dependent on $n$. }}
%with 
%\begin{equation} {C_1}^2<\frac4{\sqrt{n}}.\label{condition on C_1}\end{equation}

This $Q$ is the \emph{Cartan tensor} of Bao, Chern and Shen \cite{chern:finsler}, or the tensor $ Q$ of \cite{andrews:anisotropic} restricted to the tangent space of the level set.     The restriction here excludes anisotropic area integrands that deviate too far from the isotropic case (such as approximations to the crystalline case).  Similar restrictions on third derivatives of $\eff$ are made in studies of the elliptic problem (see, for example, \cite{winklmann:2005}).

\noindent{\textbf{Symmetry condition: }} {\textit{ This is satisfied when
\begin{equation}
\eff(p+\phi^0)=\eff(p-\phi^0) \text{ for all $p=\sum_{i=1}^n p_i\phi^i$.}\label{symmetry condition} 
\end{equation}}}

\begin{example}
 The isotropic case clearly satisfies the second condition, and also satisfies the first with $C_1=0$.  \end{example}

We impose these conditions singly, to find two time-interior gradient estimates for periodic solutions; and jointly, to find a gradient estimate interior in both space and time.    (Note that in \cite{andrews}, such estimates are found for periodic flows without imposing either condition.)

Let $u:\bigR^n\times[0,T]\rightarrow \bigR$ be a $C^3$,  bounded $|u(x,t)|\le M$, solution
to the anisotropic mean curvature flow equation \eqref{anisotropic curvature flows}.

\begin{theorem}%[First estimate for periodic anisotropic curvature flows]
\label{First estimate for periodic anisotropic curvature flows}
Let $u$ be periodic, so that $u(x,t)=u(x+L,t)$ for some lattice $L\in\bigR^n$.  If 
$\eff$ satisfies condition \eqref{condition on D^3 F} with 
${C_1}^2< 4n^{-1/2}$
 then 
 \begin{equation*}
\eff (Du-\phi^0)\le \max \left\lbrace t^{q/2}\exp\left(\frac{Aq(|u|-2M)^2}{4t}\right), P\right\rbrace
\end{equation*}
  for $0<t\le T'$, where 
 $T'$  depends on $M$, 
and $A$, $P$ and $q>1$ depend on $\eff$. 
\end{theorem}

\begin{theorem} \label{Second estimate for periodic anisotropic curve...}
Let $u$ be periodic, so that $u(x,t)=u(x+L,t)$ for some lattice $L\in\bigR^n$.  If $\eff$ satisfies the symmetry condition \eqref{symmetry condition},
then 
 \begin{equation*}
\eff (Du-\phi^0)\le \max\left\lbrace t\exp\left(\frac{A(|u|-2M)^2}{2t}\right), S \right\rbrace
\end{equation*}
  for $0<t\le T'$, where $T'$ depends on $M$, $A$ depens on  $\eff$, and $S$ depends on $\eff$ and $n$.   
\end{theorem}

\begin{theorem}[Interior estimate for anisotropic mean curvature flow] \label{interior estimate for anisotropic curve...}
Let $n>1$.  If $\eff$ satisfies both the smallness of third derivatives condition \eqref{condition on D^3 F} with $C_1^2<2/\sqrt n$,
and the symmetry condition 
\eqref{symmetry condition},
then 
 \begin{equation*}
\eff(Du-\phi^0)\le \max\left\lbrace t^{q/2}\exp\left(\frac{Aq(|u|-2M)^2}{4t}\right)\left(R^2-2\ktwo t-|x|^2\right)^{-r},P\right\rbrace
\end{equation*}
  for $0<t\le T'$.  
Here,  $A$, $P$ and $\ktwo $ depend only on $\eff$, while $T'>0$, $q>1$ and $r>1$ depend on $\eff$ and $M$.  
\end{theorem}

In the next section, we derive some technical results on $\eff$ and its derivatives.  In the third section we prove Theorems \ref{First estimate for periodic anisotropic curvature flows} and \ref{Second estimate for periodic anisotropic curve...} and in the final section we prove Theorem \ref{interior estimate for anisotropic curve...}.  

This work was part of my PhD thesis, written under the supervision of Dr. Ben Andrews at the Australian National University. I would like to thank him for many interesting discussions and helpful suggestions.

%%%%%%%%% calculations-with-the-function-F
\section{Some results regarding the function $\eff$}
 The uniform convexity 
implies that $\left.{\eff}D^2{\eff}\right|_\normal$ is positive definite  on the tangent space of the level set $\Sigma_{\eff(\normal)}$.     The homogeneity of $\eff$ leads to the disappearance of some derivatives of $\eff$ in radial directions:  
\begin{subequations}{\begin{gather}
\left.D\eff\right|_\normal (\normal )=\eff(\normal ) \label{purely radial derivatives} \\
\left.
D^2\eff\right|_\normal (\normal ,\cdot)=\left.
D^2\eff\right|_\normal (\cdot,\normal )=0\label{no radial second derivatives}  \\
\left. D(\eff D^2\eff) \right|_\normal (\normal ,\cdot,\cdot)=0. \label{third derivatives of F zero for some radial parts}
\end{gather}}\end{subequations}
These properties make it more convenient to work, not in the space $(T\mathbf{M})^*$, but rather in the tangent space to the level set $\Sigma_{\eff(\normal)}$.
Given $\normal $ normal to $\mathbf{M}$ at  $x$, we can map $v\in V^*$ to $T_{\normal}\Sigma_{\eff(\normal)}$ by setting $\widehat{v}:=v-r(v)\normal$.    The normal is not in the tangent space itself, as $\left.D\eff\right|_\normal (\normal)\not=0$, by \eqref{purely radial derivatives}.   By choosing   $r(v)$ appropriately,  $\widehat{v}$ will be in the tangent space, with
\begin{align*}
0=\left.D\eff\right|_{\normal}(\widehat{v}) 
=\left. D\eff\right|_{\normal}(v-r(v)\normal)
=\left.D\eff\right|_{ \normal}(v)-r(v) \eff (\normal),
\end{align*}
where we have used \eqref{purely radial derivatives} in the last step.  
With $r(v)=\left.D\eff\right|_{\normal}(v) / \eff (\normal)$ we then have 
$$\widehat v = v- \frac{\left.D\eff\right|_{\normal}(v)}{ \eff (\normal)}\normal,$$
which is non-zero if $v$ is not parallel to $\normal$.  Consequently, \eqref{no radial second derivatives} implies that for non-zero $v\in(T\mathbf{M})^*$,
$\left.\eff D^2\eff\right|_\normal(v,v)=\left.\eff D^2\eff\right|_\normal(\widehat v,\widehat v) >0$,  so we consider this as a new metric on $(T\mathbf{M})^*$ and write 
$G_\normal(v,w):=\left.\eff D^2\eff\right|_\normal(v,w)$.

\begin{lemma}\label{technical lemma zero} For all $P_1\ge0$,
\begin{equation*}\label{trivial consequence of convexity2} 
\eff(p-\phi^0)\ge P_1+\eff(-\phi^0) \phantom{=}\Longrightarrow\phantom{=} \eff(p)\ge P_1.
\end{equation*}
\end{lemma}
\begin{proof} 
This is a simple consequence of convexity, as $\eff(p-\phi^0)\le \eff(p)+\eff(-\phi^0)$.
\end{proof}

\begin{lemma} \label{lemma defining A and P for anisotropic case}
Let $\lbrace\phi^0,\dots,\phi^n\rbrace$ be a basis for $V^*$.
Then for each $P>\eff(-\phi^0)$ there exists $A_P>0$ such that 
$$\left. \eff D^2 \eff\right|_{p-\phi^0}(p,p)\ge A_P $$
for all $p=\sum_{i=1}^n p_i\phi^i $ with $\eff(p-\phi^0)\ge P$. 
\end{lemma}

%%%%%%%%%  proof-of-lemma-1   

%\margincomment{Proof of Lemma \ref{lemma defining A and P for anisotropic case}.}
\begin{proof}
%\margincomment{The proof of this should be much more straight-forward than it is here...} 
Set $P_1=P-\eff (-\phi^0)>0$.  By Lemma \ref{technical lemma zero},  if
 $\eff (p-\phi^0)\ge P$ then $\eff (p)\ge P_1$.     

Define $B(p):= \left. \eff  D^2 \eff \right|_{p-\phi^0}(p,p)$.   Consider this for  a fixed member of the level set  $p\in\Sigma_{P_1}\cap\text{span}\lbrace\phi^1,\dots,\phi^n\rbrace$;
%$\text{span}\lbrace\phi^1,\phi^n\rbrace \cap \lbrace q: q-\phi^0\in \Sigma_{P} \rbrace$
as $p$ is not parallel to $p-\phi^0$, $B(p)$ is positive.

Also,   \begin{align*}
\lim_{s\rightarrow \infty} B(sp)&= \lim_{s\rightarrow 0} \left.\eff  D^2\eff \right|_{sp-\phi^0}(sp,sp)
\\&=  \lim_{s\rightarrow \infty} \left.\eff  D^2\eff \right|_{sp-\phi^0}\left(sp-(sp-\phi^0),sp-(sp-\phi^0)\right)
\\&=  \lim_{s\rightarrow \infty} \left.\eff  D^2\eff \right|_{p-\phi^0/s}(\phi^0,\phi^0)
\\& =\left.\eff  D^2\eff \right|_{p}(\phi^0,\phi^0)
\\&>0,\end{align*} 
where in the second line, the additional terms added vanish according to \eqref{no radial second derivatives},  and in the third line the scaling in $s$ is allowed as $\eff D^2\eff$ is homogeneous degree zero.
The final inequality is because $\phi^0$ is not parallel to $p$.   Since $B(sp)>0$ for $1\le s < \infty$, it follows that 
$\inf_{s\in[1,\infty)} B(sp)=A_p>0$, and taking the minimum over all $p$ in the (compact and closed) level set 
gives $$\inf_{\Sigma_{P_1}\cap\text{span}\lbrace\phi^1,\dots,\phi^n\rbrace}
A_p=:{A}_P>0.$$  
\end{proof}

\begin{lemma}
If $\eff$ satisfies the symmetry condition \eqref{symmetry condition}, then homogeneity implies that
\begin{equation}
\label{symmetry condition implies...}
\begin{split}
&\left.D\eff \right|_p(\phi^0)=0  %\label{symmetry condition implies some first derivatives zero}
\\
&\left. D^2 \eff \right|_p(\phi^0,\phi^j)=0 %\notag 
\\
&\left. D^3 \eff\right|_p(\phi^0,\phi^j,\phi^k)=0 %\label{third consequence of symmetry}
\\
&\left. D^3 \eff\right|_p(\phi^0,\phi^0,\phi^0)=0, %\label{fourth consequence of symmetry}
\end{split}\end{equation} 
for all $p= \sum_{i=1}^n p_i\phi^i$ and all $j,k\not=0$.  
\end{lemma}
\begin{proof}  The symmetry condition implies that 
$$\left.D\eff\right|_p(\phi^0) =\lim_{s\rightarrow 0} \frac1s\left[\eff\left(p+{s\phi^0}/2\right)-\eff\left(p-{s\phi^0}/2\right)\right] 
=0.$$ The others may be proven similarly. 
 \end{proof}

In the following lemma, we show that the symmetry condition \eqref{symmetry condition}  can be used in a similar way to the smallness-of-third-derivatives condition \eqref{condition on D^3 F}.   We use this in the proof of Theorem \ref{Second estimate for periodic anisotropic curve...}.

\begin{lemma}
 Suppose the symmetry condition \eqref{symmetry condition} holds.   
For all $\epsilon>0$ we can find $S_\epsilon$  such that if $p=\sum_{i=1}^np_i\phi^i$ satisfies
$\eff(p-\phi^0)\ge S_\epsilon$ then  
\begin{equation} \label{some equation in lemma 9.7} 
\left| {\left.\eff D(\eff D^2\eff)\right|_{p-\phi^0}(p,\widehat q,\widehat q)}\right| 
\le \epsilon \left\lbrace\left.\eff D^2\eff\right|_{p-\phi^0}(p,p)\right\rbrace^{1/2}\left.\eff D^2\eff\right|_{p-\phi^0}(q,q)
\end{equation} 
 for all $q=\sum_{i=1}^nq_i\phi^i$.
  \label{lemma showing that symmetry is as good as small third derivatives}
\end{lemma}

%%%%%%  proof-of-lemma-3

\begin{proof} 
Let $\epsilon>0$ be given.    
As \eqref{some equation in lemma 9.7} is unchanged under the mapping $q\mapsto sq$, we consider only those $q$ on a fixed level set $\Sigma_1$.     

Our approach is to restrict $p$ to a level set  $\Sigma_P$, for $P>\eff(-\phi^0)$, and then show that under the mapping $p\mapsto sp$, the quotient
\begin{equation} \label{main eqn from lemma 3} 
\left| \frac{\left.\eff  D(\eff D^2\eff )\right|_{sp-\phi^0}(sp,\widehat q,\widehat q)}
{G(sp,sp)^{1/2}G(q,q)} \right| 
\end{equation} 
is less than $\epsilon$ for large enough $s$.  
In the above expression,
 \begin{equation*} 
\widehat q= q- \frac{\left. D \eff  \right|_{sp-\phi^0}(q)}{\eff ({sp-\phi^0})}(sp-\phi^0)=q-sr(q)p + r(q)\phi^0, \end{equation*}
with $\lim_{s\rightarrow \infty} r(q) =0$ and $\lim_{s\rightarrow\infty} sr(q)= \left.D\eff \right|_p(q)/\eff (p)$.  

As $s\rightarrow \infty$, the numerator of \eqref{main eqn from lemma 3}  is 
\begin{align*}
\lim_{s\rightarrow \infty}&  \left.\eff  D(\eff D^2\eff )\right|_{sp-\phi^0}(sp,\widehat q,\widehat q) 
 \\
&= \lim_{s\rightarrow \infty}  {\left.\eff  D(\eff D^2\eff )\right|_{p-\phi^0/s}(\phi^0,q-sr(q)p + r(q)\phi^0,q-sr(q)p + r(q)\phi^0)} \\
&={\left.\eff  D(\eff D^2\eff )\right|_{p}\left(\phi^0,q-\frac{\left.D\eff \right|_p(q)}{\eff (p)}p,q-\frac{\left.D\eff \right|_p(q)}{\eff (p)}p\right)}
\\&=0, 
\end{align*}       
by   \eqref{symmetry condition implies...}.  

If $q$ is not parallel to $p$, then the denominator is strictly positive: 
\begin{equation*}
\lim_{s\rightarrow\infty}{G_{sp-\phi^0}(sp,sp)^{1/2}G_{sp-\phi^0}(q,q)}\ge \sqrt{A_P}G_p(q,q)>0,
\end{equation*}
where $A_P$ is the constant given by Lemma \ref{lemma defining A and P for anisotropic case}.  
It follows that for each such $q$, we can find an $S$ so that $\eff (sp)\ge S$ implies  \eqref{main eqn from lemma 3} is less than $\epsilon$. 

In the case that $q$ is parallel to $p$, without loss of generality we can set $q=p$.  Multiply both numerator and denominator of \eqref{main eqn from lemma 3} by $s^2$, so that the latter is bounded below,  $\lim_{s\rightarrow\infty}G_{sp-\phi^0}(sp,sp)^{3/2} \linebreak[1] \ge {A_P}^{3/2}>0$.  The numerator is then 
\begin{equation*}
\left.\eff  D(\eff D^2\eff )\right|_{sp-\phi^0}\left(sp,[s-s^2r(p)]p+sr(p)\phi^0,[s-s^2r(p)]p+sr(p)\phi^0\right)  
\end{equation*}       
and since 
\begin{align*}
\lim_{s\rightarrow\infty}s-s^2r(p)&=\lim_{s\rightarrow\infty} s\left[1-\frac{\left. D \eff  \right|_{sp-\phi^0}(p)}{\eff ({p-\phi^0/s})}\right]\\
&=\lim_{{s'}\rightarrow 0}\frac1{s'}\left[\frac{\left. D \eff  \right|_{p}(p)}{\eff ({p})}-\frac{\left. D \eff  \right|_{p-{s'}\phi^0}(p)}{\eff ({p-{s'}\phi^0})}\right]\\
&=D\left(\frac{\left. D \eff  \right|_{p}}{\eff ({p})}\right)\left(\phi^0,p\right) \\
&= \frac{\left. D^2 \eff  \right|_{p}(\phi^0,p)}{\eff (p)} - \frac{\left. D \eff  \right|_{p}(p)\left. D \eff  \right|_{p}(\phi^0)}{\eff (p)^2}\\
&=0,
\end{align*}
the numerator approaches $\left.\eff  D(\eff D^2\eff )\right|_{p}\left(\phi^0,\phi^0,\phi^0\right)=0$ as $s\rightarrow \infty$. That is, for $q=p$ we can find an $S$ such that $\eff (sp)\ge S$ implies  \eqref{main eqn from lemma 3} is less that $\epsilon$. 
                                                                                             
Since \eqref{main eqn from lemma 3} is continuous in $q$, the supremum over $q\in\Sigma_1$ (and hence in $V^*$) of these constants $S$ is finite, and we set this to be ${S'}_\epsilon$. 
Finally, we set $S_\epsilon=S'_\epsilon+\eff (-\phi^0)$, so that whenever $\eff (p-\phi^0)\ge S_\epsilon$, $\eff (p)\ge S'_\epsilon$.  
\end{proof}

The final two technical lemmas are used in the proof of Theorem \ref{interior estimate for anisotropic curve...}.
\begin{lemma} \label{lemma giving lower bound on trace G}
Let $\lbrace \phi^0,\phi^1,\dots,\phi^n\rbrace$ be a basis for $V^*$, where $n>1$.  Then there  are constants $\kone ,\ktwo >0$ such that for all $p=\sum_{i=1}^n p_i\phi^i$,  \begin{equation*}
\kone \le \sum_{i=1}^n \left.\eff D^2\eff\right|_{p-\phi^0}(\phi^i,\phi^i)\le \ktwo .
\end{equation*}
\end{lemma}

\begin{proof}  
Let $p\in\Sigma_1$ be fixed.  Make an orthogonal change of coordinates on $\lbrace{\phi^1,\dots,\phi^n}\rbrace$ so that $p$ is parallel to $\phi^1$.  Note that $\left.G\right|_{p-\phi^0}(\phi^i,\phi^i)$ is strictly positive.   Mapping $p\mapsto sp$, we notice that
\begin{align*}&\lim_{s\rightarrow \infty}\left.G\right|_{sp-\phi^0}(\phi^i,\phi^i) \\
&\phantom{=}=  \lim_{s\rightarrow \infty}\left.G\right|_{sp-\phi^0}\left(\phi^i-\frac{\left.D\eff\right|_{sp-\phi^0}(\phi^i)}{\eff (sp-\phi^0)}(sp-\phi^0),  \phi^i-\frac{\left.D\eff\right|_{sp-\phi^0}(\phi^i)}{\eff (sp-\phi^0)}(sp-\phi^0)\right) \\
&\phantom{=}=  \left.G\right|_{p}\left(\phi^i-\frac{\left.D\eff\right|_{p}(\phi^i)}{\eff (p)}p,  \phi^i-\frac{\left.D\eff\right|_{p}(\phi^i)}{\eff (p)}p\right)\\
&\phantom{=}=\left.G\right|_{\phi^1}\left(\phi^i,  \phi^i\right),  
\end{align*}  is strictly positive and finite for $i=2,\dots,n$, as is  
\begin{align*}\lim_{s\rightarrow 0}\left.G\right|_{sp-\phi^0}(\phi^i,\phi^i)
&\phantom{=}= \left.G\right|_{-\phi^0}\left(\phi^i+\frac{\left.D\eff\right|_{-\phi^0}(\phi^i)}{\eff (-\phi^0)}\phi^0,\phi^i+\frac{\left.D\eff\right|_{-\phi^0}(\phi^i)}{\eff (-\phi^0)}\phi^0\right)\\
&\phantom{=}=\left.G\right|_{-\phi^0}\left(\phi^i,\phi^i\right).
\end{align*}

It follows that $\left.G\right|_{sp-\phi^0}(\phi^1,\phi^1) + \sum_{i=2}^n \left.G\right|_{sp-\phi^0}(\phi^i,\phi^i)$  has strictly positive bounds for all $s\in[0,\infty)$, and 
taking the minimum and maximum of these bounds over $p\in\Sigma_1$  gives the result.
\end{proof}

\begin{lemma} If $\eff$ satisfies the symmetry condition \eqref{symmetry condition},
then there exists a constant $C_2$ depending only on $\eff$ such that 
\begin{equation*}\left.\eff D^2\eff\right|_{p-\phi^0}(p,q)\le C_2\frac {\eff(q)}{\eff(p-\phi^0)}\end{equation*}
for all $p=\sum_{i=1}^n p_i\phi^i$ and  $q=\sum_{i=1}^n q_i\phi^i$.     \label{cross-terms lemma}
\end{lemma}
\begin{proof}   This  is unchanged under the mapping $q\mapsto sq$, so  we may restrict $q$ to $\Sigma_1$.  

Let $q$ be fixed.  For any given $p\in\Sigma_1$, consider
$$\frac{\eff (p-\phi^0)}{\eff (q)} \left.\eff D^2\eff\right|_{p-\phi^0}(p,q) $$
under the mapping $p\mapsto sp$ as $s\rightarrow\infty$:
\begin{align*}
\lim_{s\rightarrow\infty} &
\frac{\eff (sp-\phi^0)}{\eff (q)} 
 \left.\eff D^2\eff\right|_{sp-\phi^0}(sp,q)
\\&=\lim_{s'\rightarrow 0}\frac1{\eff (q)}\frac1{s'}\left[{\eff (p-s'\phi^0)} \left.\eff D^2\eff\right|_{p-s'\phi^0}(\phi^0,q)
-\left.\eff^2 D^2\eff\right|_{p}(\phi^0,q)\right]\\
&=-\frac1{\eff(q)}\left.D\left(\eff^2D^2\eff\right)\right|_p(\phi^0,\phi^0,q),
\end{align*}
which is bounded, as $\eff$ is $C^3$.  Consequently, 
\begin{equation*}
\sup_{s\in[0,\infty)} \frac{\eff (sp-\phi^0)}{\eff (q)} \left.\eff D^2\eff \right|_{sp-\phi^0}(sp,q)\le C(p,q) <\infty
\end{equation*}
for some finite $C(p,q)$.  Setting $C_2=\max_{p,q\in \Sigma_1} C(p,q)$ gives the required result.  
\end{proof}

%%%%% {the-gradient-estimate}

\section{The gradient estimate for periodic flows}

\begin{proof}[Proof  of Theorem \ref{First estimate for periodic anisotropic curvature flows}]
Define $Z:= \eff (Du-\phi^0)-\varphi(u,t)$
 where $\varphi$ is a smooth positive function for $t>0$ with $\varphi(\cdot,0)\ge \sup_{t=0} \eff (Du-\phi^0)$.  Later, we will choose $\varphi$ to be some inverse power of the fundamental solution to a heat equation \eqref{definition of phi}, but we start by focussing on the first part of $Z$.  
 
Consider the first point where  $Z$ is no longer negative, so that  $\eff =\varphi$.    This point will be a spatial maximum of $Z$,  due to the periodicity of $u$.

Assume that at this point, $\eff (Du-\phi^0)\ge P>\eff (-\phi^0)$.

The first derivative condition at this point is
$0=D_kZ=\left.D\eff \right.|_\normal(\phi^m)u_{mk}-{\varphi'} u_k$,
where $\normal=Du-\phi^0$.  That is, for all vectors $v\in \text{span}\lbrace e_1,\dots,e_n\rbrace$, 
\begin{equation} \label{first derivative condition for bar F}  
D^2u\left(\left.D\eff \right|_\normal(\phi^m)e_m,v\right)={\varphi'}Du\left(v\right).
\end{equation} 

Using \eqref{no radial second derivatives}, we can rewrite the evolution equation for $u$ in terms of the tangential (to the level set $\Sigma_{\eff (\normal)}$) directions $\widehat\phi^i$,
\begin{equation*}
u_t=\left.\eff  D^2\eff  \right|_\normal(\phi^i,\phi^j)u_{ij}=\left.\eff  D^2\eff  \right|_\normal(\widehat\phi^i,\widehat\phi^j)u_{ij}.
\end{equation*}

We make use of this in finding an evolution equation for $\eff $:
{\allowdisplaybreaks[1]{\begin{align}
\pd {\eff } t&= \left.D\eff \right|_\normal(\phi^k)u_{kt} \notag \\*
&= \left.D\eff \right|_\normal(\phi^k)\left[\left.\eff  D^2\eff  \right|_\normal(\widehat\phi^i,\widehat\phi^j)u_{ij}\right]_k \notag \\
&= 
\left.D\eff \right|_\normal(\phi^k)
\left[  
\left.D\left(\eff  D^2\eff \right) \right|_\normal(D_k \normal,\widehat\phi^i,\widehat\phi^j)u_{ij} +\left.\eff  D^2\eff  \right|_\normal(D_k\widehat\phi^i,\widehat\phi^j)u_{ij}
\right. \notag \\*
&\phantom{spacespacespa} \left.
+\left.\eff  D^2\eff  \right|_\normal(\widehat\phi^i,D_k\widehat\phi^j)u_{ij}   
+\left.\eff  D^2\eff  \right|_\normal(\widehat\phi^i,\widehat\phi^j)u_{ijk}
 \right] 
\notag \\* &\phantom{==}
+\left.\eff  D^2\eff \right|_{\normal}(\phi^i,\phi^j)D_{ij}F 
\notag\\* &\phantom{==}
-\left.\eff  D^2\eff \right|_{\normal}(\phi^i,\phi^j)\left[ \left.D^2\eff \right|_\normal(\phi^m,\phi^l)u_{mi}u_{lj}+\left.D\eff \right|_\normal(\phi^m)u_{mij}\right]
\notag \\
&=
\left.D\eff \right|_\normal(\phi^k)
\Big[  
\left.D\left(\eff  D^2\eff \right) \right|_\normal(D_k \normal,\widehat\phi^i,\widehat\phi^j)u_{ij} +\left.\eff  D^2\eff  \right|_\normal(D_k\widehat\phi^i,\widehat\phi^j)u_{ij}
 \notag \\*
&\phantom{spacespacespa}
+\left.\eff  D^2\eff  \right|_\normal(\widehat\phi^i,D_k\widehat\phi^j)u_{ij}  \Big]
\notag \\* &\phantom{==}
+\left.\eff  D^2\eff \right|_{\normal}(\phi^i,\phi^j)D_{ij}F 
-\left.\eff  D^2\eff \right|_{\normal}(\phi^i,\phi^j) \left.D^2\eff \right|_\normal(\phi^m,\phi^l)u_{mi}u_{lj}
\notag, 
\end{align}}} where in the third step we have added and subtracted second derivatives of $F$.
Derivatives of $\normal$ are 
$D_k\normal=u_{mk}\phi^m \label{derivative of radial z} $ 
which we use to simplify those terms with derivatives of $\widehat\phi^i$:
 \begin{align*}
\left.D^2\eff \right|_\normal(D_k\widehat\phi^i,\widehat\phi^j)&=\left.D^2\eff \right|_\normal(D_k(-
c({\phi^i})\normal),\widehat\phi^j)\\
&=\left.D^2\eff \right|_\normal(-D_k\left(c({\phi^i})\right)\normal-c({\phi^i})u_{mk}\phi^m,\widehat\phi^j) \\
&=-c({\phi^i})\left.D^2\eff \right|_\normal(u_{mk}\phi^m,\widehat\phi^j)\\
&=-\frac{\left.D\eff \right|_{\normal}(\phi^i)}{ \eff  (\normal)}\left.D^2\eff \right|_\normal(u_{mk}\phi^m,\widehat\phi^j).
\end{align*}

The evolution equation is now
\begin{align}\pd {\eff } t&=
\left.D\eff \right|_\normal(\phi^k)  
\left.D\left(\eff  D^2\eff \right) \right|_\normal(\phi^m,\widehat\phi^i,\widehat\phi^j)u_{mk}u_{ij} 
\notag 
\\ &\phantom{=}
-\left.D\eff \right|_\normal(\phi^k) \left[ \left.D\eff \right|_\normal(\phi^i)  \left. D^2\eff  \right|_\normal(\phi^m,\widehat\phi^j) +  \left.D\eff \right|_\normal(\phi^j)  \left. D^2\eff  \right|_\normal(\phi^m,\widehat\phi^i) \right]u_{mk}u_{ij}
\notag 
\\ &\phantom{=}
+\left.\eff  D^2\eff \right|_{\normal}(\phi^i,\phi^j)D_{ij}F 
-\left.\eff  D^2\eff \right|_{\normal}(\phi^i,\phi^j) \left.D^2\eff \right|_\normal(\phi^m,\phi^l)u_{mi}u_{lj}.
 \label{evolution eqn for F, first time}  
\end{align}

At a critical point of $Z$,  we can use the first derivative condition \eqref{first derivative condition for bar F} to simplify further.  The first term of  \eqref{evolution eqn for F, first time} becomes
\begin{align*}
D^2u( \left.D\eff \right|_\normal(\phi^k) e_k,e_m) 
\left.D\left(\eff  D^2\eff \right) \right|_\normal&(\phi^m,\widehat\phi^i,\widehat\phi^j)u_{ij}  \\&=
{\varphi'}Du(e_m)
\left.D\left(\eff  D^2\eff \right) \right|_\normal(\phi^m,\widehat\phi^i,\widehat\phi^j)u_{ij} \\
&= {\varphi'}
\left.D\left(\eff  D^2\eff \right) \right|_\normal(Du,\widehat\phi^i,\widehat\phi^j)u_{ij}, 
\end{align*}
while the second becomes 
\begin{align*}
&- \left.D\eff \right|_\normal(\phi^k){\left.D\eff \right|_{ \normal}(\phi^i)}\left.D^2\eff \right|_\normal(\phi^m,\widehat\phi^j)u_{mk}u_{ij}\\
&\phantom{spacespace}=- \left. D^2\eff  \right|_\normal(\phi^m,\widehat\phi^j)D^2u\left( \left.D\eff \right|_\normal(\phi^k)e_k,e_m\right)D^2u\left( \left.D\eff \right|_\normal\left(\phi^i\right)e_i,e_j\right) \\
&\phantom{spacespace}=-{\varphi'}^2\left. D^2\eff  \right|_\normal(Du(e_m)\phi^m,Du(e_j)\widehat\phi^j)\\
&\phantom{spacespace}=-{\varphi'}^2\left. D^2\eff  \right|_\normal(Du,Du),
\end{align*}
as does the third, so the evolution equation is 
\begin{align*}
\pd {\eff } t&= \frac{{\varphi'}}\varphi\left.\eff  D\left(\eff  D^2\eff \right) \right|_\normal(Du,\widehat\phi^i,\widehat\phi^j)u_{ij}-2\frac{{\varphi'}^2}\varphi  \left. \eff  D^2\eff  \right|_\normal(Du,Du)  
\notag \\ &\phantom{==}
+\left.\eff  D^2\eff \right|_{\normal}(\phi^i,\phi^j)D_{ij}F 
-\frac1\varphi\left.\eff  D^2\eff \right|_{\normal}(\phi^i,\phi^j) \left.\eff  D^2\eff \right|_\normal(\phi^m,\phi^l)u_{mi}u_{lj},
\end{align*}
where we have multiplied some terms through by $1=\eff /\varphi$ 
in order that derivatives of $\eff $ appear as homogeneous degree zero terms.   

Derivatives of $\varphi$ are given by
{\allowdisplaybreaks{\begin{gather*}
D\varphi={\varphi'} Du \\
D_{ij}\varphi={\varphi''}u_i u_j+{\varphi'} u_{ij}\\
\dfrac{d \varphi}{dt}={\varphi'}u_t+\varphi_t
\end{gather*}   }}
for $i,j\not=0$, 
so  an evolution equation for $\varphi$ is
\begin{align*}\dfrac{d \varphi}{dt}
&={\varphi'}u_t+\varphi_t +\left.\eff  D^2\eff \right|_\normal(\phi^i,\phi^j)\left( D_{ij}\varphi-{\varphi''}u_i u_j-{\varphi'} u_{ij}\right) \\
&=\varphi_t +\left.\eff  D^2\eff \right|_\normal(\phi^i,\phi^j)\left( D_{ij}\varphi-{\varphi''}u_i u_j\right),\end{align*}
and the entire evolution equation for $Z$, at a local maximum, is  
\begin{align}
\dfrac {dZ} {dt}
&=\left.\eff  D^2\eff \right|_{\normal}(\phi^i,\phi^j)D_{ij}Z+  \frac{{\varphi'}}{\varphi} \left.\eff  D\left(\eff  D^2\eff \right) \right|_{\normal}(Du,\widehat\phi^i,\widehat\phi^j)u_{ij}
\notag\\&\phantom{=} 
-2\frac{{\varphi'}^2}{\varphi}\left. \eff  D^2\eff  \right|_{ \normal}(Du,Du)  
-\frac1\varphi\left.\eff  D^2\eff \right|_{\normal}(\phi^i,\phi^j)\left.\eff  D^2\eff \right|_{\normal}(\phi^m,\phi^l)u_{mi}u_{lj} \notag
\\&\phantom{=}
-\varphi_t
+\left.\eff  D^2\eff \right|_{\normal}(Du,Du){\varphi''} .\notag
%\label{last line of evolution equation for Z, one}.
\end{align}
Notice that all the covectors $\phi^i$, $Du$  appear in places where replacing them by their projections in the  tangent space of $\Sigma_{\eff (\normal)}$,  that is, by $\widehat\phi^i$ or  $\widehat{Du}$, has no effect, thanks to  \eqref{no radial second derivatives}  and \eqref{third derivatives of F zero for some radial parts}. 
On the tangent space, $D^2\eff $ is positive definite.
Choose the basis $\lbrace \phi^1,\dots,\phi^n\rbrace$ so that $G_\normal$ is the identity at the maximum point,  ${G_\normal}^{\alpha\beta}=\delta^{\alpha \beta}$.    The evolution equation for $Z$ is now
\begin{align}
\dfrac {dZ} {dt}&=G^{ij}D_{ij}Z+\frac{{\varphi'}}{\varphi} \left.\eff  D\left(\eff  D^2\eff \right) \right|_{\normal}(\widehat{Du},\widehat\phi^i,\widehat\phi^j)u_{ij}-2\frac{{\varphi'}^2}{\varphi}G({Du},{Du})  
%&\phantom{space} 
\notag \\
&\phantom{spacespacespacespace}-\frac1\varphi G^{ij}G^{ml}u_{mi}u_{lj} -\varphi_t
 %\\ &\phantom{spacespacespace}
+G({Du},{Du}){\varphi''}.
\label{last line of evolution equation for Z}\end{align}

The Cauchy-Schwarz inequality for a positive definite matrix $B$ implies that  $v^T w \le \epsilon v^T B v +({4\epsilon})^{-1}w^T B ^{-1} w$.
We use this to estimate the second term of~\nopagebreak~\eqref{last line of evolution equation for Z}: 
%\label{page for C-s}
\begin{align*}
 \frac{{\varphi'}}{\varphi} &\left.\eff  D\left(\eff  D^2\eff \right) \right|_{ \normal}(\widehat{Du},\widehat\phi^i,\widehat\phi^j)u_{ij} \notag
%\\&= \frac{{\varphi'}}{\varphi}u_k \left.D\left(\eff  D^2\eff \right) \right|_{ z}(\widehat\phi^k,\widehat\phi^i,\widehat\phi^j)u_{ij} \\
\\&\phantom{space}= \frac{{\varphi'}}{\varphi} \left[
\left.D\eff  \right|_{ \normal} (\widehat{Du})\left. \eff  D^2\eff  \right|_{\normal}(\widehat\phi^i,\widehat\phi^j)+ \left. \eff ^2D^3\eff  \right|_{\normal}(\widehat{Du},\widehat\phi^i,\widehat\phi^j)
\right]u_{ij} \notag \\
&\phantom{space}= \frac{{\varphi'}}{\varphi}u_k {Q}^{kij}
u_{ij}  \notag \\
&\phantom{space}\le  
\epsilon \frac{{{\varphi'}}^2}\varphi G(Du,Du)
+\frac1{4\epsilon\varphi} G_{\alpha\beta }Q^{\alpha i j}u_{ij}Q^{\beta kl} u_{kl} ,  \label{cauchy-schwarz for interior estimate}
\end{align*}  
where the first term of the second line is zero, as $\widehat{Du}$ is tangent to the unit ball, so $\left.D\eff  \right|_{ \normal} (\widehat{Du})=0$.  In the last line, we have used the notation for the inverse $G_{\alpha\beta}=(G^{-1})^{\alpha\beta}$.  

We can use \eqref{condition on D^3 F}, the smallness-of-third-derivatives condition, to estimate the second term in this inequality: 
\allowdisplaybreaks[1]{\begin{align*}
\frac1{4\epsilon\varphi} G_{\alpha\beta} Q^{\alpha i j}u_{ij}Q^{\beta kl} u_{kl} 
&=\frac1{4\epsilon\varphi} Q\left(G_{\alpha\beta}\widehat\phi^\alpha,u_{ij}\widehat\phi^i,\widehat\phi^j\right)
 Q\left(\widehat\phi^\beta,u_{kl}\widehat\phi^k,\widehat\phi^l\right)\\*
&\le\frac{ {C_1}^2}{4\epsilon\varphi} 
\bigg( G(G_{\alpha \beta}\widehat\phi^\alpha,G_{\gamma \beta} \widehat\phi^\gamma)G(u_{ij}\widehat\phi^i,u_{mj}\widehat\phi^m)G(\widehat\phi^j,\widehat\phi^j) %\right. 
\\*&\phantom{spacespace}%\left.
\times G(\widehat\phi^\beta, \widehat\phi^\beta)G(u_{kl}\widehat\phi^k,u_{pl}\widehat\phi^p)G(\widehat\phi^l,\widehat\phi^l)
\bigg)^{1/2} \\
&=
\frac{ {C_1}^2}{4\epsilon\varphi} 
\left( G_{\beta\beta}u_{ij}u_{ij}G^{jj}G^{\beta\beta}u_{kl}u_{kl}G^{ll}\right)^{1/2}\\*
&=\frac{ {C_1}^2}{4\epsilon\varphi} 
\sqrt{n}  \left( G^{ij}G^{kl}u_{ik}u_{jl}\right).
\end{align*}}

Now we can estimate \eqref{last line of evolution equation for Z}  from above ---
\begin{align*}
\dfrac {dZ} {dt}&\le 
G^{ij}D_{ij}Z
+\frac1\varphi\left(\frac{ {C_1}^2}{4\epsilon} 
\sqrt{n}    -1\right) G^{ij}G^{kl}u_{ik}u_{jl}\\&\phantom{spacespacespace} 
+\frac{{{\varphi'}}^2}\varphi \left(\epsilon-2\right)G(Du,Du) 
-\varphi_t\notag +{\varphi''}G(Du,Du).
\end{align*}

The second term is zero if we choose 
$\epsilon\,=\,{C_1}^2\sqrt{n}/4\,<\,1$.

Choose $\varphi=\Phi^{-q}$ for some $q>1$ and
\begin{equation} \label{definition of phi}\Phi(u,t)=\frac1{\sqrt t}\exp\left(-A\frac{(u-2M)^2}{4t}\right),\end{equation}
which satisfies the heat equation $\Phi_t=A\Phi''$, where $A=A_P$ is the constant given by Lemma  \ref{lemma defining A and P for anisotropic case}.   As 
\begin{gather*}
\varphi'=-q\Phi^{-q-1}\Phi'\notag \\
\varphi''=q(q+1)\Phi^{-q-2}(\Phi')^2 -q\Phi^{-q-1}\Phi'' \notag \\
\varphi_t=-q\Phi^{-q-1}\Phi_t, \notag 
\end{gather*}
the equation satisfied by $\varphi$ is
$\varphi_t= A\varphi'' -A\left(1+q^{-1}\right){\varphi'^2}/\varphi.
\notag$  
If we substitute $\Phi$ and its derivatives for $\varphi$ and its derivatives, we find that
\begin{align*}
\frac{{{\varphi'}}^2}\varphi &\left(\epsilon\right.-\left.2\right)G(Du,Du) 
-\varphi_t\notag +{\varphi''}G(Du,Du)
\\&\phantom{=} =
q^2\Phi^{-q-2}{\Phi'}^2(\epsilon-2)G(Du,Du)\\&\phantom{space}
+\left[q(q+1)\Phi^{-q-2}{\Phi'}^2-q\Phi^{-q-1}\Phi''\right]G(Du,Du)
        +Aq\Phi^{-q-1}\Phi'' \\
&\phantom{=}=q\Phi^{-q-2}{\Phi'}^2\left[q(\epsilon-1)+1\right]G(Du,Du)+ q\Phi^{-q-1}\Phi''\left[A-G(Du,Du)\right]. 
\end{align*}

The first term is zero if we choose 
$q^{-1}={1-\epsilon}={1-{C_1}^2\sqrt{n}/4}.$

As we assumed at the beginning 
  that $\eff (Du-\phi^0)\ge P$, Lemma \ref{lemma defining A and P for anisotropic case}
implies that 
$G(Du,Du)\ge A_P.$
As $\Phi''$ is positive for small times, for $t<T'$ we have $\pd Z t \le 0$.

On the other hand, if we consider the possibility that $\eff(Du-\phi^0)=\varphi<  P$ at this local maximum, we could replace $\varphi$ by $\sup \lbrace \varphi, P \rbrace$ in the definition of $Z$.    In that case, the first maximum of $Z$ occurs at a point where the barrier is flat, and so the first variation is  
$$0=D_kZ=\left.D\eff \right|_z(\phi^k)u_{mk},$$
and the evolution equation for $Z$ at the local maximum is
\begin{align*} 
\frac{dZ}{dt}&= \left.\eff  D^2\eff \right|_\normal(\phi^i,\phi^j)D_{ij}Z-\left. D^2\eff \right|_{\normal}(\phi^i,\phi^j) \left.\eff  D^2\eff \right|_\normal(\phi^m,\phi^l)u_{mi}u_{lj} 
\le 0.
\end{align*}
 
Since $Z_t\le 0$ at the first point where $Z=0$,  $Z\le 0 $ for all $t<T'$.  The same argument works if, in the definition of $\Phi$ \eqref{definition of phi}, the term $u-2M$ is replaced by $u+2M$.  The conclusion (with $|u|-2M$)  follows.  
\end{proof}

\begin{proof}[Proof of Theorem \ref{Second estimate for periodic anisotropic curve...}]
We begin by defining $Z$ as in Theorem \ref{First estimate for periodic anisotropic curvature flows}, and assume that the first non-negative value of $Z$ occurs when $\eff (Du-\phi^0)\ge S_{\epsilon}$, for some $\epsilon>0$  to be chosen later and the corresponding $S_\epsilon$ given by Lemma \ref{lemma showing that symmetry is as good as small third derivatives}.   We then follow the earlier proof up to 
equation \eqref{last line of evolution equation for Z}, the evolution equation for $Z$ at a local maximum.

This time, we choose the local coordinates $\lbrace \phi^1,\dots,\phi^n\rbrace$ so that at this point $D^2u$ is diagonal.    This puts the second term of \eqref{last line of evolution equation for Z} in a suitable form to be estimated using Lemma \ref{lemma showing that symmetry is as good as small third derivatives}.    
\begin{align*}
\frac{{\varphi'}}{\varphi} \left.\eff  D\left(\eff  D^2\eff \right) \right|_{\normal}(Du,\widehat\phi^j,\widehat\phi^j)u_{jj}
&\le \left|\epsilon \frac{{\varphi'}}{\varphi} \sqrt{G(Du,Du)}G(\widehat\phi^j,\widehat\phi^j)u_{jj}\right|
\\&\le \frac{{\varphi'}^2}{2\varphi}G(Du,Du) + \frac{{\epsilon}^2}{2\varphi}\left|G(\widehat\phi^j,\widehat\phi^j)u_{ij}\right|^2
\\&\le \frac{{\varphi'}^2}{2\varphi}G(Du,Du) + \frac{{\epsilon}^2}{2\varphi}n G^{ij}u_{jk}G^{kl}u_{li},
\end{align*}
where in the last line we use trace inequality $(\trace A)^2\le n\trace (A^2)$.  

If we now choose $\epsilon=\sqrt{2/n}$, the second term of this inequality is cancelled by the fourth term of \eqref{last line of evolution equation for Z}.   The evolution equation now becomes
\begin{align*}
\dfrac {dZ} {dt}
&\le G^{ij}D_{ij}Z-\frac32 \frac{{{\varphi'}}^2}{\varphi}G(Du,Du) 
-\varphi_t\notag +G({Du},{Du}){\varphi''} .
\end{align*}
This is negative at a local maximum if we make the same choice of barrier as before, $\varphi=\Phi^{-q}$ for $q=2$, $\Phi$ given by \eqref{definition of phi}, with $A=A_{S_\epsilon}$ given by Lemma \ref{lemma defining A and P for anisotropic case}.  

If our assumption that $\eff (Du-\phi^0)\ge S_\epsilon$ does not hold, then we can replace $\varphi$ by $\max\lbrace S_\epsilon,\varphi\rbrace$.  At the local maximum, $Z_t\le 0$ and so the conclusion follows.
\end{proof}

\begin{remark} In the last theorem, we have chosen $q=2$ somewhat arbitrarily;  in fact $q$ needs only to be strictly greater than $1$, since we can set $q= (1-n\epsilon^2/4)^{-1}$, for  $\epsilon$ given by Lemma \ref{lemma showing that symmetry is as good as small third derivatives}.  However, a smaller $\epsilon$ may force a larger $S_\epsilon$, so the optimal choice would depend on the exact form of $\eff $. 
\end{remark}

%%%%%%%%%% {the-interior-gradient-estimate}

\section{Interior estimate for anisotropic mean curvature flow}

\begin{proof}[Proof  of Theorem \ref{interior estimate for anisotropic curve...}]
We introduce a localising term $\eta$ into our definition of $Z$,
\begin{equation*}
Z:=\eff (Du-\phi^0)-\frac\varphi\eta,  
\end{equation*} 
which is now restricted to the shrinking ball $(x,t)\in B_{\sqrt{R^2-2\ktwo t}}\times[0,T]$, where $\ktwo $ is the constant given by Lemma \ref{lemma giving lower bound on trace G}.  The smooth strictly positive function $\varphi=\varphi(u,t)$ is   chosen so that $Z<0$ at the initial time, and $\eta$ is a smooth positive function chosen so that $\eta\rightarrow 0$ on the boundary of the shrinking ball.

Assume that at the first interior point where $Z=0$, $\eff (Du-\phi^0)\ge P>\eff (-\phi^0)$. 

Then $F(Du)=\varphi/\eta$ and as this is a spatial maximum (since the choice of $\eta$ ensures that there are no boundary maxima)
we have a first derivative condition
\begin{equation}  \label{first derivative condition for interior anisotropic estimate}
0=D_kZ=\left.D\eff \right.|_\normal(\phi^m)u_{mk}-D_k\left(\varphi/\eta\right).
\end{equation}

An evolution equation for  $\varphi/\eta$, with second derivatives $G^{ij}D_{ij}\left(\varphi/\eta\right)$ added and subtracted, follows:
\begin{align*}
\dfrac{d}{dt} \left(\frac \varphi \eta\right)
&=\frac1\eta\left({\varphi'} u_t+\varphi_t\right) -\frac\varphi{\eta^2}\frac{d\eta}{dt}+ \left.\eff  D^2\eff \right|_{\normal}(\phi^i,\phi^j)D_{ij}\left(\frac \varphi \eta\right)
%\\&\phantom{space}
\\&\phantom{===}
-\left.\eff  D^2\eff \right|_{\normal}(\phi^i,\phi^j)\bigg[\frac1\eta\left({\varphi''}u_iu_j+{\varphi'}u_{ij}\right)-\frac{{\varphi'}}{\eta^2}\left(u_jD^i\eta+u_iD^j\eta\right)
\\&\phantom{======================}
+2\frac\varphi{\eta^3}D^i\eta D^j\eta-\frac{\varphi}{\eta^2}D_{ij}\eta \bigg]\\
&=G^{ij}D_{ij}\left(\frac\varphi\eta\right)+\frac1\eta\left[\varphi_t-G(Du,Du){\varphi''} \right]
-\frac\varphi{\eta^2}\left(\frac{d}{dt}- G^{ij}D_{ij}\right)\eta 
\\&\phantom{===}
+2\frac{{\varphi'}}{\eta^2}G(Du,D\eta) -2\frac\varphi{\eta^3}G(D\eta,D\eta).
\end{align*}

We can incorporate the first derivative condition \eqref{first derivative condition for interior anisotropic estimate} into \eqref{evolution eqn for F, first time},  the evolution equation for $\eff $: {\allowdisplaybreaks{
\begin{align*}\dfrac {d\eff }{d t}&=
G^{ij}D_{ij}\eff   + 
D_m\left(\varphi/\eta\right)  
\left.D\left(\eff  D^2\eff \right) \right|_\normal(\phi^m,\widehat\phi^i,\widehat\phi^j)u_{ij} 
\notag 
\\*&\phantom{=i}
-2\left.\eff  D^2\eff \right|_\normal \left(D(\varphi/\eta),D(\varphi/\eta)\right)
-\left.\eff  D^2\eff \right|_{\normal}(\phi^i,\phi^j) \left.D^2\eff \right|_\normal(\phi^m,\phi^l)u_{mi}u_{lj} 
 \\
&=G^{ij}D_{ij}\eff   + \frac{\eta}{\varphi} \left.\eff 
D\left(\eff  D^2\eff \right) \right|_\normal\left(
D(\varphi/\eta),\widehat\phi^i,\widehat\phi^j\right)u_{ij} \\*
&\phantom{=i}-2\frac{\eta}{\varphi}\left[\frac{{\varphi'}^2}{\eta^2}G(Du,Du)-2\frac{\varphi{\varphi'}}{\eta^3}G(Du,D\eta)+\frac{\varphi^2}{\eta^4}G\left(D\eta,D\eta\right)\right] 
\notag \\*  &\phantom{=i}
-\frac{\eta}{\varphi}G^{ij}G^{ml}u_{mi}u_{lj}.
\end{align*} }}

Putting the last two steps together gives an evolution equation for $Z$ at a local maximum:
\begin{equation}\label{starred equation}
\begin{split}
\dfrac {dZ}{d t} &= 
%\dfrac {d\eff }{d t} -\dfrac{d }{dt}\left(\frac\varphi\eta\right) \\ 
G^{ij}D_{ij}Z +\frac{\eta}{\varphi} 
\left.\eff D\left(\eff  D^2\eff \right) \right|_\normal\left(D(\varphi/\eta),\widehat\phi^i,\widehat\phi^j\right)u_{ij} 
\\&\phantom{=} -\frac{\eta}{\varphi}G^{ij}G^{ml}u_{mi}u_{lj}
-\frac1\eta\left[\varphi_t-G(Du,Du){\varphi''} +2\frac{{\varphi'}^2}{\varphi}G(Du,Du)\right] 
\\&\phantom{=}+\frac\varphi{\eta^2}\left(\frac{d}{dt}-G^{ij}D_{ij}\right)\eta
+2\frac{{\varphi'}}{\eta^2}G(Du,D\eta).
\end{split}
\end{equation}

The second term here may be split up into a part with $D\varphi$ and a part with $D\eta$:
\begin{align*}
\frac{\eta}{\varphi} &\left.\eff 
D\left(\eff  D^2\eff \right) \right|_\normal \left(\frac{{\varphi'}}\eta Du-\frac{\varphi}{\eta^2} D\eta,\widehat\phi^i,\widehat\phi^j\right)u_{ij} \\
&=\frac{{\varphi'}}{\varphi} \left.\eff 
D\left(\eff  D^2\eff \right) \right|_\normal 
        \left(Du,\widehat\phi^i,\widehat\phi^j\right)u_{ij}  
-\frac{1}{\eta} \left.\eff  D\left(\eff  D^2\eff \right) \right|_\normal 
        \left(D\eta,\widehat\phi^i,\widehat\phi^j\right)u_{ij}. 
\end{align*}

These may be individually estimated using the Cauchy-Schwarz inequality and the smallness-of-third-derivatives condition, 
as described in the proof of Theorem \ref{First estimate for periodic anisotropic curvature flows} ---
{\allowdisplaybreaks{\begin{gather*}
\begin{split}
\frac{{\varphi'}}{\varphi} \left.\eff 
D\left(\eff  D^2\eff \right) \right|_\normal 
        &\left(Du,\widehat\phi^i,\widehat\phi^j\right)u_{ij}  \\*
&\le \mu_1 \frac{{\varphi'}^2}{\varphi\eta} G(Du,Du)  + \frac1{4\mu_1}\frac\eta{\varphi} C_1^2\sqrt{n} \left(G^{ij}G^{ml}u_{mi}u_{lj}\right),
\end{split}\\
\begin{split}
-\frac{1}{\eta} \left.\eff  D\left(\eff  D^2\eff \right) \right|_\normal& 
\left(D\eta,\widehat\phi^i,\widehat\phi^j\right)u_{ij}\\*
%&\le \frac{\mu_2\varphi} {\eta^3}G(D\eta,D\eta)+\frac\eta{4\mu\varphi}G_{\alpha \beta}Q^{\alpha i j }u_{ij}Q^{\beta k l}u_{kl}\\
&\le \mu_2\frac \varphi {\eta^3}G(D\eta,D\eta)+\frac1{4\mu_2}\frac\eta{\varphi}C_1^2\sqrt n\left(G^{ij}G^{kl}u_{ik}u_{jl}\right),
\end{split}\end{gather*} }}
for some $0<\mu_1,\mu_2<1$.

We choose the localising term $\eta:=\tilde \eta^r$ for some $r>1$ and $\tilde \eta=R^2-2\ktwo t-|x|^2$.  Then
$D_i\eta= r\tilde \eta^{r-1}D_i\tilde\eta$,%\\
    $D_{ij}\eta=r\tilde \eta^{r-1}D_{ij}\tilde \eta + r(r-1)\tilde \eta^{r-2}D_i\tilde \eta D_j\tilde \eta$,
and  the second-last term of the evolution equation \eqref{starred equation} is 
\begin{align*}
\frac\varphi{\eta^2}\left(\frac{d}{dt}-G^{ij}D_{ij}\right)\eta&=
\frac\varphi{\eta^2}r{\tilde\eta}^{r-1}\left[-2\ktwo +2\trace G -(r-1){\tilde\eta}^{-1}G(D\tilde\eta, D\tilde\eta) \right]\\
&\le \frac\varphi{\eta^2}r{\tilde\eta}^{r-2}(1-r)G(D\tilde\eta, D\tilde\eta). \end{align*}

As $\eff $ satisfies the symmetry condition \eqref{symmetry condition}, we may use  Lemma \ref{cross-terms lemma} to estimate the
 final term of the evolution equation:
\begin{align*}
2\frac{{\varphi'}}{\eta^2}G(Du,D\eta)&=2\frac{{\varphi'}}{\eta^2}\left.\eff D^2\eff \right|_{Du-\phi^0}(Du,D\eta)
\\&\le 2\frac{{\varphi'}}{\eta^2} \frac{C_2\eff (D\eta)}{\eff (Du-\phi^0)}\\
&= 2C_2\eff (D\eta)\frac{{\varphi'}}{\varphi\eta}\\
&\le 2C_2C_3rR^{2r-1}\frac{{\varphi'}}{\varphi\eta},
\end{align*}
where we have used that $\eff (D\eta)= r\tilde\eta^{r-1}\eff (D\tilde\eta)\le C_3rR^{2r-1}$, for  $C_3>0$ depending only on $\eff $.

The evolution equation can now be estimated from above: 
\begin{align}\dfrac {dZ}{d t}
&\le G^{ij}D_{ij}Z   + \frac{\eta}\varphi\left(\frac1{4\mu_1} C_1^2\sqrt n+\frac1{4\mu_2}C_1^2\sqrt n -1 \right) G^{ij}G^{ml}u_{mi}u_{lj} 
\notag \\&\phantom{space}
-\frac1\eta\left[\varphi_t-G(Du,Du){\varphi''} +(2-\mu_1)\frac{{\varphi'}^2}{\varphi}G(Du,Du)
-2C_2\eff (D\eta)\frac{{\varphi'}}\varphi\right]
\notag \\&\phantom{spacespacespace}+\frac\varphi{\eta^2}r\tilde\eta^{r-2} \left(1-r+r\mu_2\right)G(D\tilde \eta,D\tilde \eta). \label{almost the last eqn}
\end{align}

Since ${C_1^2\sqrt n}/4<1/2$, we can choose $\mu_1<1$ and $\mu_2<1$ 
such that 
\begin{gather*}
\frac{C_1^2\sqrt{n}}{4}\left(\frac1\mu_1+\frac1\mu_2\right)\le 1. %\\
\end{gather*}

With such choices, the second term of the evolution inequality \eqref{almost the last eqn} will be negative.  We can also  set $r=(1-\mu_2)^{-1}>1$,  so the coefficient of $\tilde\eta^{-1}G(D\tilde \eta,D\tilde \eta)$ is zero.

As in the previous cases we set $\varphi=\Phi^{-q}$ where $\Phi$ is given by \eqref{definition of phi} for $A=A_P$ given by Lemma \ref{lemma defining A and P for anisotropic case}.

The bracketted part  of the second line of \eqref{almost the last eqn} is then
\begin{align*} 
&-\frac1\eta\left[\varphi_t-G(Du,Du){\varphi''} +(2-\mu_1)\frac{{\varphi'}^2}{\varphi}G(Du,Du)
-2C_2\eff (D\eta)\frac{{\varphi'}}\varphi \right] \notag\\
&\phantom{=====}=\frac1\eta\Bigg[q\Phi^{-q-1}\Phi_t+ G(Du,Du)\left(q(q+1)\Phi^{-q-2}{\Phi'}^2-q\Phi^{-q-1}\Phi''\right)\\
&\phantom{========}-(2-\mu_1)G(Du,Du)q^2\Phi^{-q-2}{\Phi'}^2+2C_1C_2rR^{2r-1}\Phi^{-1}|\Phi'|\Bigg]\\ 
&\phantom{=====}=\frac{q\Phi^{-q-1}}\eta\left(\Phi_t-G(Du,Du)\Phi''\right) \\
&\phantom{========}+\frac{q|\Phi'|}{\eta\Phi}\left[G(Du,Du)|\Phi'|\Phi^{-q-1}\left(1-q+q\mu_1\right) +2C_1C_2rR^{2r-1} \right].  
\end{align*}

If we choose $T'$ small enough that $\Phi''\ge0$, then the term $\Phi_t-G(Du,Du)\Phi'' =\linebreak\left(A-G(Du,Du)\right)\Phi''$ is negative.
Additionally, if we choose $T'$ small enough that $\Phi\le1$, we have 
\begin{equation*}
G(Du,Du)|\Phi'|\Phi^{-q-1}\ge \frac{A^2}tM\Phi^{-q}\ge \frac{A^2M}{T'},\end{equation*}
and so we need only to choose $q$ large enough that 
\begin{equation*}
q\ge \frac1{1-\mu_1}\left(1+\frac{2C_1C_2rR^{2r-1}T'}{A^2M}\right)
\end{equation*}
for the last term to be negative.
So, at such maxima, $Z_t\le0$.  

At local maxima where $\eff  (Du-\phi^0)<P$, then in the definition of $Z$ we replace $\varphi/\eta$ by $\max\lbrace \varphi/\eta, P\rbrace$, in which case the barrier is flat at the local maxima, and we again find that 
$Z_t\le0$.

In either case, the maximum principle ensures that $Z$ is never greater than zero and the conclusion follows. 
\end{proof}

\end{document}